\documentclass[english,letterpaper]{amsart}
\usepackage[T1]{fontenc}
\usepackage[latin1]{inputenc}
\usepackage{verbatim}
\usepackage{amssymb}

\makeatletter

\providecommand{\tabularnewline}{\\}

 \theoremstyle{plain}    
 
 \numberwithin{equation}{section} 
 \numberwithin{figure}{section} 

\usepackage{hyperref}

\usepackage{mathrsfs}

\usepackage{mathtools}

\usepackage{euscript}

\usepackage{upgreek}

\usepackage{babel}
\makeatother
\begin{document}
\newcommand{\dual}[1]{#1^{\vee}}

\newcommand{\trans}[1]{^{t}\!#1}

\newcommand{\sym}[2]{\mathrm{Sym}^{#1}#2}

\newcommand{\alt}[2]{\mathrm{Alt}^{#1}#2}

\newcommand{\set}[2]{\left\{  #1\, |\, #2\right\}  }

\newcommand{\map}[3]{#1\!:#2\!\rightarrow\!#3}

\newcommand{\aut}[1]{\mathrm{Aut}\!\left(#1\right)}

\newcommand{\End}[1]{\mathrm{End}\!\left(#1\right)}

\newcommand{\id}[1]{\mathbf{id}_{#1}}

\newcommand{\ip}[1]{\left[#1\right]}

\newcommand{\uhp}{\mathbf{H}}

\newcommand{\sm}[4]{\left(\begin{smallmatrix}#1  &  #2\cr\cr#3  &  #4\end{smallmatrix}\right)}

\newcommand{\cyc}[1]{\mathbb{Q}\left[\zeta_{#1}\right]}

\newcommand{\zz}{\exp\left(2\pi i/3\right)}

\newcommand{\ZN}[1]{\left(\mathbb{Z}/#1\mathbb{Z}\right)^{\times}}

\newcommand{\Mod}[3]{#1\equiv#2\, \left(\mathrm{mod}\, \, #3\right)}

\newcommand{\psl}{\mathrm{PSL}_{2}\!\left(\mathbb{Z}\right)}

\newcommand{\SL}{\mathrm{SL}_{2}\!\left(\mathbb{Z}\right)}

\newcommand{\gl}[2]{\mathrm{GL}_{#1}\!\left(#2\right)}

\newcommand{\tr}[2]{\textrm{Tr}_{_{#1}}\left(#2\right)}

\newcommand{\mr}[4]{\rho\!\left(\begin{array}{cc}
 #1  &  #2\\
#3  &  #4\end{array} \right)}

\newcommand{\FA}{\mathbb{X}}

\newcommand{\FB}{\mathbf{\Lambda}}

\newcommand{\FC}{\mathfrak{z}}

\newcommand{\am}{\EuScript{A}}

\newcommand{\bm}{\EuScript{B}}

\newcommand{\dm}[1]{\mathfrak{D}\!\left(#1\right)}

\newcommand{\FD}{\mathbb{X}\left(q,z\right)}

\newcommand{\FE}{\EuScript{P}}

\newcommand{\FF}[1]{\EuScript{E}\!\left(#1\right)}

\newcommand{\mro}[1]{\mathcal{M}\!\left(#1\right)}

\newcommand{\SB}[2]{\FA^{\left(#1;#2\right)}}

\newcommand{\gfxi}{\mathfrak{X}_{\xi\eta}\left(q,z\right)}

\newcommand{\ev}[1]{{\bf e}_{#1}} \newcommand{\pcsg}[1]{\Gamma\left(#1\right)}

\newcommand{\ks}{\lambda}

\newcommand{\orb}{\eta}

\newcommand{\ccorb}{\mathcal{O}}

\newcommand{\DO}{\nabla} \newcommand{\ei}{\mathcal{E}} \newcommand{\xm}[3]{\mathcal{X}_{#1}^{\left(#2;#3\right)}}
\newcommand{\om}[1]{\mathbf{\Xi}\!\left(#1\right)}

\newcommand{\mo}[1]{\mathbf{\Xi}^{#1}\!\left(\tau\right)}

\newcommand{\D}[1]{#1^{\prime}\!\!} \newcommand{\diff}{\mathrm{\mathsf{d}}}

\newcommand{\dif}[2]{\frac{\mathsf{d}#1}{\mathsf{d}#2}} \newcommand{\sig}[1]{\upsigma\!\left(#1\right)}

\newcommand{\kan}{\varkappa}

\newcommand{\gam}[1]{\upgamma_{#1}}

\newcommand{\mof}[2]{M_{#1}\!\left(#2\right)}

\newcommand{\cof}[2]{S_{#1}\!\left(#2\right)}

\newcommand{\varmu}{\mu}

\author{P. Bantay and T. Gannon}

\email{bantay@general.elte.hu}

\address{Institute for Theoretical Physics, Eötvös Loránd University, Budapest}

\email{tgannon@math.ualberta.ca}

\address{Department of Mathematical Sciences, University of Alberta, Edmonton}


\title{Vector-valued modular functions for the modular group and the 
hypergeometric equation}

\thanks{Work of P.B. was supported by grants OTKA T047041, T043582, the János
Bolyai Research Scholarship of the Hungarian Academy of Sciences and
EC Marie Curie MRTN-CT-2004-512194. T.G. would like to thank E\"otv\"os University and the University of Hamburg
for kind hospitality while this research was undertaken; his research is
supported in part by NSERC and the Humboldt Foundation.}

\begin{abstract}
A general theory of vector-valued modular functions, holomorphic in the 
upper half-plane, is presented for
finite dimensional representations of the modular group. This also
provides a description of vector-valued modular forms of arbitrary half-integer
weight. It is shown that the space of these modular functions is 
spanned, as a module over the polynomials in $J$, 
by the columns of a matrix that satisfies an abstract
hypergeometric equation, providing a simple solution of the Riemann-Hilbert
problem for representations of the modular group. Restrictions on
the coefficients of this differential equation implied by analyticity are
discussed, and an inversion formula is presented that allows the determination
of an arbitrary vector-valued modular function from its singular behavior.
Questions of rationality and positivity of expansion coefficients
are addressed. Closed expressions for the number of vector-valued
modular forms of half-integer weight are given, and the general theory
is illustrated on simple examples.
\end{abstract}
\maketitle

\section{Introduction}

The notions of modular functions and forms -- and their generalizations
-- are among the most fruitful in all of mathematics, and with the
arrival of String Theory they have become standard fare in mathematical
physics as well. Vector-valued modular functions $\FA(\tau)$ for
$\SL$ appear for instance as characters of Vertex Operator Algebras
\cite{bib:Zh} and Conformal Field Theories \cite{bib:CFT}, and in
the Norton series of generalized Moonshine \cite{bib:No1};
moreover in Conformal Field Theory,
 vector-valued modular forms of arbitrary rational weight appear
as conformal blocks on a once-punctured torus. In spite of its importance,
there has been little attempt at a systematic treatment of this theory
(\cite{bib:KM,bib:ES} are exceptions).

In these contexts, singularities of the component functions $\FA_{\orb}(\tau)$
appear at the cusps $\mathbb{Q}\cup\{\infty\}$, but not in
the upper half-plane $\uhp$, and we will restrict our attention to
such functions. In a previous paper \cite{bib:BG} we explained (with
examples) how to obtain all such vector-valued modular functions,
given the corresponding multiplier $\rho$, a finite-dimensional representation
of $\left(\mathrm{P}\right)\SL$. In this paper we focus on the underlying
structure of these spaces of vector-valued modular functions.
They are generated by the $\SL$-Hauptmodul $J\!\left(\tau\right)$,
together with the columns of a certain \textit{fundamental matrix}
$\om{\tau}$.
We explain how everything is conveniently recovered from the exponents
$\FB$ at infinity and a numerical matrix $\mathcal{X}$ (essentially,
the first nontrivial $q$-coefficients of $\om{\tau}$). The other
$q$-coefficients of $\om{\tau}$ can be obtained from a differential
equation, the monodromy of which is determined by $\rho$. Our results
extend directly to vector-valued modular forms of half-integer weight:
for instance, we obtain an explicit formula for the dimension of
the spaces of such forms.

In Section \ref{sec:The-fundamental-matrix}, we review the framework
of \cite{bib:BG}, and discuss a subtlety: the choice of integer part
of the exponent matrix $\FB$.
 Section \ref{sec:The-hypergeometric-form} explains how the differential
equation satisfied by the fundamental matrix may be recast into an
abstract hypergeometric equation, and the consequences this has on
the various quantities involved. Section \ref{sec:Low-dimensional-examples}
gives some concrete examples, illustrating the effectiveness of our
results. Section \ref{sec:The-inversion-formula} provides an inversion
formula, which allows the explicit computation of any vector-valued
modular function 
from its singular part, provided the fundamental matrix is known.
In the motivating examples, the $q$-expansions have nonnegative integer
coefficients: Section \ref{sec:Positivity-and-integrality} explains
how the existence of such $q$-expansions constrains $\rho$. An appendix
describes what happens when -- as is typical in Vertex Operator Algebras
or Conformal Field Theory -- $\rho$ is a representation of $\SL$
rather than of $\mathrm{P}\SL$.

\section{The fundamental matrix\label{sec:The-fundamental-matrix}}

Consider a matrix representation $\map{\rho}{\SL}{\gl{d}{\mathbb{C}}}$
whose kernel contains $\sm{-1}{\;\:0}{\;\,0}{-1}$, and for which
$T=\rho\sm{1}{1}{0}{1}$ is a diagonal matrix of finite order. 
We associate to $\rho$ the set $\mro{\rho}$ of all those maps $\map{\FA}{\uhp}{\mathbb{C}^{d}}$
which are holomorphic in the upper half-plane $\uhp=\set{\tau}{\mathrm{Im}\tau>0}$,
transform according to $\rho$, that is %
\footnote{Here and in what follows we view $\FA\left(\tau\right)$ as a column
vector.%
} \begin{equation}
\FA\left(\frac{a\tau+b}{c\tau+d}\right)=\mr{a}{b}{c}{d}\FA\left(\tau\right)\label{eq:modtrans}\end{equation}
 for all $\sm{a}{b}{c}{d}\in\SL$ and $\tau\in\uhp$, and have only
finite order poles at the cusps \cite{bib:BG}. This last condition
means the following: since $\rho\sm{1}{1}{0}{1}$ is diagonal of finite
order, there exists a diagonal matrix $\FB$ (the \textit{exponent
matrix}) such that \begin{equation}
\mr{1}{1}{0}{1}=\exp\left(2\pi i\FB\right)\:,\label{eq:expdef}\end{equation}
 the diagonal elements of $\FB$ being rational numbers. Because of
Eq.\eqref{eq:modtrans}, the map $\exp\left(-2\pi i\tau\FB\right)\FA\left(\tau\right)$
is periodic in $\tau$ (with period 1): consequently, it may be expanded
into a Fourier series%
\footnote{In all what follows, we shall alternate freely between the notations
$f\left(\tau\right)$ and $f\left(q\right)$ for one and the same
quantity $f$: in general, the notation $f\left(\tau\right)$ is meant
to emphasize that we consider $f$ as a (holomorphic) function on
the upper half-plane $\uhp$, while $f\left(q\right)$ refers to its
expansion as a power series in $q=\exp\left(2\pi i\tau\right)$. %
} \begin{equation}
q^{-\FB}\FA\left(\tau\right)=\sum_{n\in\mathbb{Z}}\FA\left[n\right]q^{n}\:,\label{eq:qexp}\end{equation}
 where $q=\exp\left(2\pi i\tau\right)$. We define the principal part
$\FE\FA$ of $\FA$ as the sum of the terms with negative powers of
$q$ on the rhs. of Eq.\eqref{eq:qexp}, i.e. \begin{equation}
\FE\FA\left(q\right)=\sum_{n<0}\FA\left[n\right]q^{n}\:.\label{eq:Pdef}\end{equation}
With this definition, $\FA$ has finite order poles at the cusps if
and only if its principal part $\FE\FA$ is a finite sum.

Clearly, the space $\mro{\rho}$ is an infinite dimensional linear
space over $\mathbb{C}$, a basis being provided by the maps $\SB{\xi}{n}\in\mro{\rho}$
which have a pole of order $n>0$ at the $\xi$th position, i.e. \begin{equation}
\left[\FE\SB{\xi}{n}\left(q\right)\right]_{\eta}=q^{-n}\delta_{\xi\eta}\:.\label{eq:sbdef}\end{equation}
 We call these $\SB{\xi}{n}$ the \textit{canonical basis vectors};
they are clearly linearly independent, and that they exist and therefore
span $\mro{\rho}$ was explained in \cite{bib:BG} (an independent
proof is provided at the end of Section \ref{sec:The-hypergeometric-form}).

Let\begin{equation}
J\!\left(\tau\right)=q^{-1}+\sum_{n=1}^{\infty}c\left(n\right)q^{n}=q^{-1}+196884q+\cdots\label{eq:Jexp}\end{equation}
 denote the Hauptmodul of $\SL$, i.e. the (suitably normalized) generator
of the field of modular functions for $\SL$ (for this and other aspects
of the classical theory of modular functions and forms, see e.g. \cite{bib:Ap}).
Multiplication by $J$ takes the space $\mro{\rho}$ to itself, in
other words $\mro{\rho}$ is a $\mathbb{C}\left[J\right]$-module.
The important point is that this is a (free) $\mathbb{C}\left[J\right]$-module
of \textit{finite} rank, because the canonical basis vectors satisfy
the \textit{recursion relations} \cite{bib:BG}
\begin{equation}
\SB{\xi}{m+1}=J\!\left(\tau\right)\SB{\xi}{m}-\sum_{n=1}^{m-1}c\left(n\right)\SB{\xi}{m-n}-\sum_{\orb}\xm{\eta}{\xi}{m}\SB{\orb}{1}\,\,,\label{eq:recursion}\end{equation}
 where \begin{equation}
\xm{\eta}{\xi}{m}=\SB{\xi}{m}\left[0\right]_{\eta}=\lim_{q\rightarrow0}\left(\left[q^{-\FB}\SB{\xi}{m}\left(q\right)\right]_{\eta}-q^{-m}\delta_{\xi\eta}\right)\label{eq:xmdef}\end{equation}
 denotes the {}``constant part'' of $\SB{\xi}{m}$. These recursion
relations allow to express each canonical basis vector $\SB{\xi}{m}$
in terms of the $\SB{\xi}{1}$-s, proving that the latter generate
the $\mathbb{C}\left[J\right]$-module $\mro{\rho}$. Later on, we'll
give an explicit expression -- Eq.\eqref{eq:cbasis} -- for the $\SB{\xi}{m}$-s.
We will see shortly that the $\SB{\xi}{1}$ are linearly independent
over the field $\mathbb{C}\left(J\right)$ of modular functions, and
thus the $\mathbb{C}\left[J\right]$-module $\mro{\rho}$ has rank
$d$.

Besides the recursion relations Eq.\eqref{eq:recursion}, there is
a second set of relations -- the \textit{differential relations} \cite{bib:BG}
-- between the canonical basis vectors. They follow from the fact
that the differential operator \begin{equation}
\DO=\frac{\FF{\tau}}{2\pi i}\dif{}{\tau}\:\label{eq:nabladef}\end{equation}
 maps $\mro{\rho}$ to itself, where\begin{equation}
\FF{\tau}=\frac{E_{10}\left(\tau\right)}{\Delta\left(\tau\right)}=\sum_{n=-1}^{\infty}\mathcal{E}_{n}q^{n}=q^{-1}-240-141444q-\cdots\,\:\label{eq:Hdef}\end{equation}
 is the quotient of the (normalized) Eisenstein series of weight 10
by the discriminant form $\Delta\left(\tau\right)=q\prod_{n=1}^{\infty}\left(1-q^{n}\right)^{24}$
of weight 12. Looking at the action of $\DO$ on the canonical basis
vectors, one gets the differential relations \begin{equation}
\DO\SB{\xi}{m}=\left(\FB_{\xi\xi}-m\right)\sum_{n=-1}^{m-1}\ei_{n}\SB{\xi}{m-n}+\sum_{\orb}\FB_{\eta\orb}\xm{\eta}{\xi}{m}\SB{\orb}{1}\,\,.\label{eq:difrel}\end{equation}

The compatibility of the recursion and differential relations requires
that %
\footnote{One may show, using the results of Section \ref{sec:The-inversion-formula},
that this is not only a necessary, but also a sufficient condition
for the compatibility of the recursion and differential relations.%
}

\begin{equation}
\DO\SB{\xi}{1}=\left(J-240\right)\left(\FB_{\xi\xi}-1\right)\SB{\xi}{1}+\sum_{\eta}\left(1+\FB_{\eta\eta}-\FB_{\xi\xi}\right)\xm{\eta}{\xi}{1}\SB{\eta}{1}\:,\label{eq:compat1}\end{equation}
 which is a first order ordinary differential equation -- the \textit{compatibility
equation} -- for the $\SB{\xi}{1}$-s.

One may recast the compatibility equation Eq.\eqref{eq:compat1} in
a more suggestive form by introducing the \textit{fundamental matrix}
\begin{equation}
\om{\tau}_{\xi\eta}=\left[\SB{\eta}{1}\left(\tau\right)\right]_{\xi}\:,\label{eq:funddef}\end{equation}
 whose columns span over $\mathbb{C}\left[J\right]$ the module $\mro{\rho}$.
Then Eq.\eqref{eq:compat1} takes the form \begin{equation}
\frac{1}{2\pi i}\dif{\om{\tau}}{\tau}=\om{\tau}\dm{\tau}\:,\label{eq:compat}\end{equation}
 where \begin{equation}
\dm{\tau}=\frac{1}{\FF{\tau}}\left\{ \left(J\!\left(\tau\right)-240\right)\left(\FB-1\right)+\mathcal{X}+\left[\FB,\mathcal{X}\right]\right\} \;\label{eq:Ddef}\end{equation}
 and $\mathcal{X}_{\xi\eta}=\xm{\xi}{\eta}{1}$ is the so-called \textit{characteristic
matrix} (as usual, $\left[\mathcal{X},\FB\right]=\mathcal{X}\FB-\FB\mathcal{X}$
denotes the commutator of matrices). Note that Eq.\eqref{eq:compat}
has singular points at the poles of $\dm{\tau}$, i.e. at the $\SL$-orbits
of the cusp $\tau=i\infty$ and elliptic points $\tau=i$ and $\tau=\zz$.
Taking into account the boundary condition \begin{equation}
q^{1-\FB_{\xi\xi}}\om{q}_{\xi\eta}=\delta_{\xi\eta}+O\left(q\right)\;\mathrm{as}\;\mbox{ }q\rightarrow0\:,\label{eq:xibc}\end{equation}
which follows from Eq.\eqref{eq:sbdef}, one can solve Eq.\eqref{eq:compat},
provided one knows the exponent matrix $\FB$ and the characteristic
matrix $\mathcal{X}$, determining then from Eq.\eqref{eq:recursion}
the canonical basis vectors $\SB{\xi}{m}$. The theory of ordinary
differential equations guarantees Eq.\eqref{eq:compat} to have series
solutions that converge in suitably small neighborhoods of $\uhp$
avoiding the elliptic points, but the holomorphicity of $\om{\tau}$
implies that those series actually converge throughout $\uhp$.

Eq.\eqref{eq:xibc} tells us that the determinant det$\,\om{\tau}$
has leading term $q^{\tr{}{\FB-1}}$ as $q\rightarrow0$, and so is
not identically 0. Thus, its columns $\SB{\xi}{1}$ are indeed linearly
independent over $\mathbb{C}\left(J\right)$. This invertibility of
$\om{\tau}$ legitimates its appellation, since it is now seen as
a fundamental solution of Eq.\eqref{eq:compat}.

Actually, the results so far enable us already to discuss vector-valued
modular forms of half-integer weight for $\SL$. By a modular form
of weight $k\in\frac{1}{2}\mathbb{Z}$ for the (possibly projective)
P$\SL$-representation $\varrho$ we'll mean a map $\FA:\uhp\rightarrow\mathbb{C}^{d}$
that is holomorphic everywhere in $\uhp$, transforms according to
\begin{equation}
\FA\left(\frac{a\tau+b}{c\tau+d}\right)=\left(c\tau+d\right)^{k}\varrho
\!\left(\begin{array}{cc} a  &  b\\ c  &  d\end{array} \right)
\FA\left(\tau\right)\:,\label{eq:modformtrans}\end{equation}
 and which tends to a finite limit as $\tau\rightarrow i\infty$.
Such an $\FA$ is a cusp form if it vanishes at $\tau=i\infty$. As
before, we require $\varrho\sm{1}{1}{0}{1}$ to be diagonal and of
finite order. We'll denote by $\mof{k}{\varrho}$ and $\cof{k}{\varrho}$
the space of vector-valued modular forms (resp. cusp forms) of weight
$k$ for the representation $\varrho$: clearly, the latter is a subspace
of the former. Note that when $\varrho$ is the trivial representation, we
recover the classical theory of modular forms of even weight.

Let $\eta(\tau)=q^{1/24}\prod_{n=1}^{\infty}\left(1-q^{n}\right)$
be the Dedekind eta function, and let $\varmu$ denote its multiplier
 (see e.g.\ Chapter 4 of \cite{bib:Knop} for a formula for $\varmu$).
Then, for any $k\in\frac{1}{2}\mathbb{Z}$  we have a natural injective
map\begin{equation}
\begin{aligned}\phi:\mof{k}{\varrho} & \rightarrow\mro{\varrho\otimes\varmu^{-2k}}\\
\FA & \mapsto\eta^{-2k}\FA\;\;,\end{aligned}
\label{eq:phidef}\end{equation}
 which allows to embed the space $\mof{k}{\varrho}$ of modular forms
of weight $k$ into the module $\mro{\varrho\otimes\varmu^{-2k}}$.
It is easy to verify that the space $\mof{k}{\varrho}$ will be trivial
unless $\rho=\varrho\otimes\varmu^{-2k}$ is a true representation
of P$\SL$, in which case $\mro{\rho}$ is covered by the analysis
of this paper. Nonsingularity of $\FA\in\mof{k}{\varrho}$ as $\tau\rightarrow i\infty$
bounds the order of the pole of the singular part of $\eta^{-2k}\FA$;
as a result, the spaces $\mof{k}{\varrho}$ (hence $\cof{k}{\varrho}$
too) are finite dimensional, and explicit bases can be found. 
As for their dimension, one obtains the result \begin{equation}
\begin{aligned}\dim\mof{k}{\varrho}= & \max\left(0,\mathrm{Tr}\ip{\FB+\frac{k}{12}}\right)\:,\\
\dim\cof{k}{\varrho}= & \max\left(0,-\mathrm{Tr}\ip{1-\frac{k}{12}-\FB}\right)\:,\end{aligned}
\label{eq:dimformula}\end{equation}
where $\FB$ denotes the exponent matrix of $\rho=\varrho\otimes\varmu^{-2k}$,
and $\left[x\right]$ denotes the integer part of $x$ ($x$ can be
a matrix): note that $\FB$ varies  with the weight $k$.

When $\varrho$ is the trivial representation, Eq.\eqref{eq:dimformula}
reduces to classical results for the dimensions of modular and cusp
forms for $\SL$. Those equations also lead to the following expressions
for the trace of the integer part of $\FB$ (for a true P$\SL$ representation
$\varrho$): \begin{equation}
\mathrm{Tr}\ip{1-\FB}=\dim\mof{2}{\overline{\varrho}}\label{eq:trL1}\end{equation}
 and \begin{equation}
\mathrm{Tr}\ip{\FB}=\dim\mof{0}{\varrho}-\dim\cof{2}{\overline{\varrho}}\:,\label{eq:trL}\end{equation}
 where $\overline{\varrho}$ denotes the contragredient representation
of $\varrho$. We leave the derivation of these results -- which amount
to careful bookkeeping -- to a future publication. Eq.\eqref{eq:dimformula}
recovers and generalizes the dimension formula in \cite{bib:ES},
which was proved using the Eichler-Selberg trace formula.

At this point we should make an important proviso: Eq.\eqref{eq:expdef}
only determines the fractional part of the diagonal elements of the
exponent matrix, not their integer part. This is important, since
the values of these integer parts enter the definition Eq.\eqref{eq:Pdef}
of the principal part map $\FE$, hence of the canonical basis vectors
$\SB{\xi}{n}$. Another choice of these integer parts leads to a different
set of canonical basis vectors, hence different characteristic and
fundamental matrices, while $\mro{\rho}$ remains unchanged. Even
more important is the observation that for an arbitrary choice of
the integer part of $\FB$, the principal part map $\FE$ may not
be injective (i.e. the terms singular with respect to $\FB$ may not
determine the functions) and may not be surjective (i.e. not all canonical
basis vectors may exist). As we are going to explain, one can choose
the integer part of the exponent matrix at will, provided that the
relation%
\footnote{Here we assume that the matrix representation $\rho$ is indecomposable,
i.e. cannot be written as the direct sum of two matrix representations
(this holds for any representation coming from e.g. RCFT): otherwise,
one should apply these considerations to each direct summand separately. %
} \begin{equation}
\tr{}{\FB}=\frac{5d}{12}+\frac{1}{4}\tr{}{S}+\frac{2}{3\sqrt{3}}\mathrm{Re}\left(e^{-\frac{\pi i}{6}}\tr{}{U}\right)\:\label{eq:Riemann-Roch}\end{equation}
 holds, where $d$ is the dimension of $\rho$, and we use the notations
$S=\rho\!\sm{0}{-1}{1}{0}$ and $U=\rho\!\sm{0}{-1}{1}{-1}$. Surjectivity
fails if Tr$(\FB)$ is greater than the rhs. of Eq.\eqref{eq:Riemann-Roch};
injectivity fails if it is less.

To see how this comes about, first note that $\FE$ is invertible iff both 
$\mathbf{\Xi}$ exists and its columns span $\mro{\rho}$.
Suppose that $\mathcal{X}_{12}\neq0$,
and consider the matrix \begin{equation}
M\left(\tau\right)=\left(\begin{array}{ccccc}
0 & -\mathcal{X}_{12} & 0 & \cdots & 0\\
\frac{1}{\mathcal{X}_{12}} & J\!\left(\tau\right)-C & -\frac{\mathcal{X}_{13}}{\mathcal{X}_{12}} & \cdots & -\frac{\mathcal{X}_{1d}}{\mathcal{X}_{12}}\\
0 & -\mathcal{X}_{32} & 1 & 0 & 0\\
\vdots & \vdots & 0 & \ddots & 0\\
0 & -\mathcal{X}_{d2} & 0 & 0 & 1\end{array}\right)\:,\label{eq:trafoM}\end{equation}
 where $C$ is a suitable constant. All matrix elements of $M$ belong
to $\mathbb{C}\left[J\right]$, and the same holds for the inverse
matrix $M^{-1}$, since $\det M=1$ irrespectively of the value of
the constant $C$. Consequently, the columns of the matrix $\mathbf{\Xi}^{\prime}\!\left(\tau\right)=\om{\tau}M\left(\tau\right)$
span $\mro{\rho}$ over $\mathbb{C}[J]$, iff those of $\mathbf{\Xi}$ do.
 By a suitable choice of
the constant $C$ one can achieve that $\mathbf{\Xi}^{\prime}\!\left(\tau\right)$
satisfies the boundary condition Eq.\eqref{eq:xibc} with \begin{equation}
\FB^{\prime}=\FB+\left(\begin{array}{ccccc}
1\\
 & -1\\
 &  & 0\\
 &  &  & \ddots\\
 &  &  &  & 0\end{array}\right)\:.\label{eq:trafoL}\end{equation}
 This means that $\mathbf{\Xi}^{\prime}\!\left(\tau\right)$ is also
a fundamental matrix for $\mro{\rho}$, corresponding to the exponent
matrix $\FB^{\prime}$. 

By applying a suitable sequence of transformations of the above type,
we see that one can add to $\FB$ any integral diagonal matrix whose
trace vanishes. But can we alter the trace of $\FB$ as well? The
answer is no, for we'll see in Section \ref{sec:The-hypergeometric-form}
(when we'll have all the necessary tools at our disposal) that the
invertibility of $\FE$ implies Eq.\eqref{eq:Riemann-Roch}.

In summary, the structure of the $\mathbb{C}\left[J\right]$-module
$\mro{\rho}$ is completely determined by the fundamental matrix $\om{\tau}$,
once an exponent matrix $\FB$ satisfying Eqs.\eqref{eq:expdef} and
\eqref{eq:Riemann-Roch} has been chosen. The fundamental matrix is
itself completely determined by the pair $\left(\FB,\mathcal{X}\right)$
of exponent and characteristic matrices, namely as the solution of
the compatibility equation Eq.\eqref{eq:compat} satisfying the boundary
condition Eq.\eqref{eq:xibc}. For this reason, we consider the pair
$\left(\FB,\mathcal{X}\right)$ as the basic data characterizing the
representation $\rho$.

For example, the representation $\rho$ may be recovered from the
compatibility equation. Indeed, Eq.\eqref{eq:compat} is invariant
under modular transformations\[
\tau\mapsto\frac{a\tau+b}{c\tau+d}\]
 for $\sm{a}{b}{c}{d}\!\in\!\SL$, which means that such a transformation
takes a solution to another solution. Since the equation is linear,
this new solution is of the form $M\om{\tau}$ for some matrix $M\!\in\!\gl{d}{\mathbb{C}}$.
Comparing this with Eqs.\eqref{eq:modtrans} and \eqref{eq:funddef},
and using the aforementioned invertibility of $\om{\tau}$, we see
that $M=\rho\sm{a}{b}{c}{d}$, the matrix representing the given modular
transformation, i.e. \begin{equation}
\om{\frac{a\tau+b}{c\tau+d}}=\rho\sm{a}{b}{c}{d}\om{\tau}\:.\label{eq:xitrans}\end{equation}
 But this argument works for any pair $\left(\FB,\mathcal{X}\right)$,
i.e. any such pair determines a finite dimensional representation
$\rho$ via Eq.\eqref{eq:xitrans}. This seems to suggest that the
pair $\left(\FB,\mathcal{X}\right)$ could be chosen at will, but
this is not the case: the analyticity of the fundamental matrix $\om{\tau}$
-- namely, that it is single valued and holomorphic throughout the
whole upper half-plane $\uhp$, including the elliptic points $\tau=i$
and $\tau=\zz$ -- puts severe restrictions on the pair $\left(\FB,\mathcal{X}\right)$.
To understand these, it turns out to be convenient to transform the
compatibility equation to an equivalent form, which is the subject
of the next section.

\section{The hypergeometric form of the compatibility equation\label{sec:The-hypergeometric-form}}

Consider the function \begin{equation}
\FC\!\left(\tau\right)=\frac{984-J\!\left(\tau\right)}{1728}\,\,,\label{eq:zdef}\end{equation}
 which maps the upper half-plane $\uhp$ onto the complex plane $\mathbb{C}$.
Note that $\FC\!\left(i\right)=0$ and $\FC\!\left(e^{2\pi i/3}\right)=1$.
As usual, we extend the definition of $\FC$ so that it maps $\tau=i\infty$
to $\infty$ (since $\FC$ has a first order pole at the cusp $\tau=i\infty$).
$\FC$ is clearly modular invariant, i.e. it maps points on the same
$\SL$ orbit to the same point of $\mathbb{C}$, and can thus be viewed
as a map from $\uhp/\SL$ to $\mathbb{C}$. Viewed this way, it is
one-to-one, and at the elliptic points $\tau=i$ and $\tau=\zz$ it
has valence 2 (respectively 3) -- this smooths the conical singularities
of the modular curve $\uhp/\SL$. Finally, $\FC\!\left(\tau\right)$
satisfies the differential equation\begin{equation}
\DO\FC=1728\:\FC\left(\FC-1\right)\:.\label{eq:zdif}\end{equation}
 The simplest way to see that Eq.\eqref{eq:zdif} holds is to note
that $\DO\FC$ is modular invariant, holomorphic in $\uhp$, and has
a pole of order 2 at $\tau=i\infty$, hence it is a quadratic polynomial
in $\FC$; moreover, it vanishes at the elliptic points because $E_{10}$
vanishes there. Eq.\eqref{eq:zdif} then follows by comparing the
coefficients of $q^{-2}$.

Let's consider the fundamental matrix as a (multivalued) function
of $\FC\:$. Then, by applying the chain rule and Eq.\eqref{eq:zdif},
one arrives at the following form of the compatibility equation:\begin{equation}
\dif{\om{\FC}}{\FC}=\om{\FC}\left(\frac{\am}{2\FC}+\frac{\bm}{3\left(\FC-1\right)}\right)\,\,,\label{eq:hyper}\end{equation}
 with \begin{subequations}\begin{eqnarray}
\am & = & \frac{31}{36}\left(1-\FB\right)-\frac{1}{864}\left(\mathcal{X}+\left[\FB,\mathcal{X}\right]\right)\label{eq:ABdefa}\\
\bm & = & \frac{41}{24}\left(1-\FB\right)+\frac{1}{576}\left(\mathcal{X}+\left[\FB,\mathcal{X}\right]\right)\,\,.\label{eq:ABdefb}\end{eqnarray}
 \end{subequations}

The important observation is that Eq.\eqref{eq:hyper} is an abstract
hypergeometric equation, since it has three regular singular points
(at $\FC=0,1$ and $\infty$), and much is known about the analytic
properties of the solutions of Eq.\eqref{eq:hyper} (background for
the following material is provided in e.g. Chapter 6 of \cite{bib:Hi}).
As a function of $\FC$ the fundamental matrix is not single valued
-- its multivaluedness, i.e. the monodromy of Eq.\eqref{eq:hyper},
is described by the representation $\rho$. In particular, the mo\-no\-dromies
around $\FC=0,\FC=1,\FC=\infty$ are given by $S=\rho\!\sm{0}{-1}{1}{0}$,
$U=\rho\!\sm{0}{-1}{1}{-1}$, $T=\rho\!\sm{1}{1}{0}{1}$ respectively.
Because the residues of Eq.\eqref{eq:hyper} at these points are $\am/2$,
$\bm/3$ and $\FB-1$, the matrices $S$ and $U$ are conjugate to
$\exp\left(\pi i\am\right)$ and $\exp\left(2\pi i\bm/3\right)$ res\-pectively,
and one has $SU=T^{-1}=\exp\left(-2\pi i\FB\right)$. We find that
the monodromy group of the abstract hypergeometric equation Eq.\eqref{eq:hyper}
is precisely the image of $\rho$.

Let's concentrate on the singular points $\FC=0$ and $\FC=1$ of
Eq.\eqref{eq:hyper}. The denominators 2 and 3 of the residues $\am/2$
and $\bm/3$ match the valence of the corresponding elliptic points.
Since the fundamental matrix is by definition holomorphic in the whole
upper half-plane, in particular at the elliptic points, it follows
that

\begin{enumerate}
\item the matrices $\am$ and $\bm$ are simple (i.e. diagonalizable), since
otherwise $\om{\tau}$ would have logarithmic singularities; 
\item the eigenvalues of $\am$ and $\bm$ are nonnegative to avoid poles; 
\item the eigenvalues of $\am$ and $\bm$ are integers, otherwise $\om{\tau}$
would have (algebraic or transcendental) branch points; 
\item no two eigenvalues of $\am/2$ and $\bm/3$ may differ by nonzero
integers, otherwise one would get logarithmic branch points. 
\end{enumerate}
These already restrict the matrices $\am$ and $\bm$ to a great extent,
but there is one more restriction, namely that all eigenvalues of
$\am/2$ and $\bm/3$ should be less than 1. This last condition is
a completeness condition: would there be an eigenvalue greater or
equal to one, the columns of the solution of Eq.\eqref{eq:hyper}
would not span the full $\mathbb{C}\left[J\right]$-module $\mro{\rho}$
(for the monodromy representation $\rho$). More precisely, let $P^{-1}AP$
be a diagonal matrix $D$, and suppose $D_{\eta\eta}\ge2$; then the
$\eta$-th column of $\om{\FC}P$ will be a multiple of $\FC$. This
column vector, as a function of $\tau$, could be divided by $J\!\left(\tau\right)-984$
while remaining holomorphic; the quotient would still be in $\mro{\rho}$,
but not in the $\mathbb{C}\left[J\right]$-span of the columns of
$\om{\tau}$. The argument for $\bm/3$ is similar, using $J\!\left(\tau\right)+744$
instead.

This last completeness condition, together with the four analyticity
conditions, imply the following

\textbf{Spectral condition}: \textsl{the possible eigenvalues of $\am$
are 0 or 1, while those of $\bm$ are either 0, 1 or 2}.

This is a basic result, which restricts considerably the possible
coefficient matrices. In particular, it implies that the characteristic
polynomials of $\am$ and $\bm$ read

\begin{equation}
\begin{aligned}\det\left(z-\am\right)=\: & z^{d-\alpha}\left(z-1\right)^{\alpha}\:,\\
\det\left(z-\bm\right)=\: & z^{d-\beta_{1}-\beta_{2}}\left(z-1\right)^{\beta_{1}}\left(z-2\right)^{\beta_{2}}\:,\end{aligned}
\label{eq:ABdet}\end{equation}
 where $d$ denotes their dimension, while the multiplicities $\alpha,\,\beta_{1}$
and $\beta_{2}$ are given by \begin{align}
\alpha=\,\, & \tr{}{\am}\:,\nonumber \\
\beta_{1}=\,\, & 2\tr{}{\bm}-\tr{}{\bm^{2}}\:,\label{eq:sigdef}\\
\beta_{2}=\,\, & \frac{1}{2}\left(\tr{}{\bm^{2}}-\tr{}{\bm}\right)\:.\nonumber \end{align}
 The quadruple $\left(d,\alpha,\beta_{1},\beta_{2}\right)$ of nonnegative
integers is a very important discrete invariant of the representation
$\rho$, which we'll call its \textit{signature}. For example, the
traces of the representation matrices $S=\rho\!\sm{0}{-1}{1}{0}$
and $U=\rho\!\sm{0}{-1}{1}{-1}$ are completely determined by it%
\footnote{Conversely, the traces of $S$ and $U$ -- together with the dimension
$d$ -- determine the signature.%
}:\begin{align}
\tr{}{S} & =\, d-2\alpha\:,\label{eq:traceS}\\
\tr{}{U} & =\, d-\frac{3}{2}\left(\beta_{1}+\beta_{2}\right)+i\frac{\sqrt{3}}{2}\left(\beta_{1}-\beta_{2}\right)\:.\label{eq:traceST}\end{align}
 We also note that \begin{align}
\tr{}{\mathcal{X}} & =4\left(62\beta_{1}+124\beta_{2}-123\alpha\right),\quad\label{eq:traceX}\\
\tr{}{\FB} & =d-\frac{\alpha}{2}-\frac{\beta_{1}+2\beta_{2}}{3}\:.\label{eq:tracelambda}\end{align}
 In particular, the trace of the characteristic matrix $\mathcal{X}$
is always an integer divisible by $4$, which is congruent to $4\alpha$
modulo $248$.

As another application of the notion of signature, let's mention the
following formula for the determinant of the fundamental matrix: \begin{equation}
\det\om{\tau}=\left(\frac{E_{4}\left(\tau\right)}{\Delta\left(\tau\right)^{1/3}}\right)^{\beta_{1}+2\beta_{2}}\left(\frac{E_{6}\left(\tau\right)}{\Delta\left(\tau\right)^{1/2}}\right)^{\alpha}\:,\label{eq:detxi}\end{equation}
 where $E_{4}$ and $E_{6}$ denote the (normalized) Eisenstein series
of weights 4 and 6. The proof of this result is simple: since $\om{\tau}$
satisfies Eq.\eqref{eq:compat}, its determinant satisfies -- according
to a theorem of Liouville -- 
the differential equation \begin{equation}
\frac{1}{2\pi i}\frac{\diff\left(\log\det\om{\tau}\right)}{\diff\tau}=\mathrm{Tr}\,\dm{\tau}\,\,.\label{eq:detdif}\end{equation}
 Moreover, it follows from Eq.\eqref{eq:xibc} that $\det\om{q}$
behaves as $q^{\tr{}{\FB-1}}$ for $q\rightarrow0$. It is an easy
matter to check that the rhs. of Eq.\eqref{eq:detxi} satisfies the
differential equation Eq.\eqref{eq:detdif} with this particular boundary
condition, and by general theory such a solution is unique.

It follows in particular that the fundamental matrix is invertible
everywhere except the elliptic points. That it can't be invertible
at the elliptic points, for typical representations,
is obvious: for example, at $\tau=i$ one has $\om{i}=S\om{i}$ because
of Eq.\eqref{eq:xitrans}, so $\om{i}$ invertible would imply $S$
trivial.

Let's return to the spectral condition. It follows from Eq.\eqref{eq:ABdet} that the minimal polynomials
of $\am$ and $\bm$ divide $z\left(z-1\right)$, resp. $z\left(z-1\right)\left(z-2\right)$.
Since any matrix is a root of its minimal polynomial, the spectral
condition may be expressed as \begin{equation}
\am\left(\am-1\right)=\bm\left(\bm-1\right)\left(\bm-2\right)=0\:.\label{eq:minipoly}\end{equation}

Of the four matrices $\FB,\,\mathcal{X},\,\am$ and $\bm$, any two
determine the other two%
\footnote{This is trivial unless two eigenvalues of $\FB$ differ by $1$, but
this can be always avoided by the use of transformations as in Eq.\eqref{eq:trafoM}. %
}, e.g. Eqs.(\ref{eq:ABdefa},b) imply that $\bm=3\left(1-\FB-\am/2\right)$.
Inserting this expression into Eq.\eqref{eq:minipoly}, one gets the
following system of algebraic equations:\begin{equation}
\begin{split}\am^{2} & =\am\:,\\
\am\FB\am=-\frac{17}{18}\am-2\left(\am\FB^{2}+\FB\am\FB+\FB^{2}\am\right) & +3\left(\am\FB+\FB\am\right)-4\FB^{3}+8\FB^{2}-\frac{44}{9}\FB+\frac{8}{9}\;.\end{split}
\label{eq:monodromy}\end{equation}
 That is, for a given exponent matrix $\FB$, the matrix $\am$ has
to satisfy Eq.(\ref{eq:monodromy}): note that this is a simultaneous
system of quadratic equations for the matrix elements of $\am$, and
that the matrix $\FB$ (which plays the role of a parameter) is diagonal.
Once a solution to Eq.\eqref{eq:monodromy} is known, the corresponding
characteristic matrix may be determined from Eq.\eqref{eq:ABdefa}.

What can be said about the solutions of Eq.\eqref{eq:monodromy}?
First of all, if $\left(\FB,\mathcal{X}\right)$ is a solution and
$M$ is a monomial matrix (i.e. the product of a diagonal and a permutation
matrix), then $\left(M^{-1}\FB M,M^{-1}\mathcal{X}M\right)$ is again
a solution: more generally, this holds for any matrix $M$, provided
that $M^{-1}\FB M$ is still diagonal. These transformations do not
change the equivalence class of the corresponding representation $\rho$,
and may be used to put the solution into some useful standard form.

More interesting is duality, the involutive transformation $\left(\FB,\mathcal{X}\right)\mapsto\left(\dual{\FB},\dual{\mathcal{X}}\right)$
with%
\footnote{We denote by $\trans{M}$ the transpose of a matrix $M$.%
} \begin{align}
\dual{\FB}=\; & \frac{5}{6}-\FB\label{eq:duallambda}\\
\dual{\mathcal{X}}=\; & 4-\,\trans{\mathcal{X}}\:,\label{eq:dualX}\end{align}
 which sends $\am$ to $\dual{\am}=1-\,\trans{\am}$ and $\bm$ to
$\dual{\bm}=2-\,\trans{\,\bm}$: clearly, $\dual{\am}$ and $\dual{\bm}$
satisfy the spectral condition if $\am$ and $\bm$ did. The fundamental
matrix corresponding to the dual pair $\left(\dual{\FB},\dual{\mathcal{X}}\right)$
is given by \begin{equation}
\mo{\vee}=\frac{E_{14}\left(\tau\right)}{\Delta^{7/6}\left(\tau\right)}\left(^{t}\om{\tau}\right)^{-1}\:.\label{eq:dualxi}\end{equation}
The prefactor is needed to ensure holomorphicity, which can be proved
using Eq.\eqref{eq:detxi} and the spectral condition. The dual representation
$\dual{\rho}$ is equivalent to the tensor product of the contragredient
of $\rho$ with the 1-dimensional representation $\kan$ appearing
in the bottom row of Table 1 below.

It is now time to establish 
the relation of invertibility of $\FE$ to Eq.\eqref{eq:Riemann-Roch},
left pending in Section \ref{sec:The-fundamental-matrix}. 
 If $\FE$ is invertible, then a fundamental matrix $\om{\tau}$ satisfying
Eq.\eqref{eq:xibc} exists for which the whole theory presented above
holds. Comparing Eqs.(\ref{eq:traceS}),(\ref{eq:traceST}) and Eq.\eqref{eq:tracelambda},
we arrive at Eq.\eqref{eq:Riemann-Roch}. In other words, while the
integer part of $\FB$ is to a great extent arbitrary, its trace is
completely determined by the representation $\rho$.

More generally, given any $\mathbb{C}[J]$-submodule $M$ of $\mro{\rho}$
of full rank $d$, linear algebra shows how to construct a matrix
$\om{\tau}$ of form Eq.\eqref{eq:funddef}, for some choice of $\FB$,
such that $M$ is the $\mathbb{C}[J]$-span of the columns of $\om{\tau}$.
Moreover, Tr$(\FB)$ will be bounded above by the rhs. of Eq.\eqref{eq:Riemann-Roch},
with strict inequality if $M\ne\mro{\rho}$ (to see this, use transformations
like Eq.\eqref{eq:trafoM} to make the $\FB$-s for $\mro{\rho}$
and $M$ agree in all but one spot). If in addition the submodule
is $\DO$-stable, then that matrix $\om{\tau}$ will satisfy Eq.\eqref{eq:hyper}
for $\am,\bm$ defined by Eqs.(\ref{eq:ABdefa},b), although the eigenvalues
of $\am$ and $\bm$ can now be arbitrary nonnegative integers. However,
this submodule can be {}``completed'' using the method outlined
in our proof of the spectral condition, by dividing the appropriate
vectors by $J-984$ or $J+744$ (at each stage, the submodule will
be $\DO$-stable, thanks to Eqs.\eqref{eq:zdif},\eqref{eq:hyper}).
The result will be matrices $\om{\tau},\FB,\am,\bm$ satisfying the
spectral condition and Eq.\eqref{eq:Riemann-Roch}. To summarize,
given a $\DO$-stable rank $d$ submodule $M$ of $\mro{\rho}$, with
matrices $\FB,\am,\bm$, we have: $M=\mro{\rho}$ iff $\FB$ satisfies
Eq.\eqref{eq:Riemann-Roch}, iff $\am,\bm$ satisfy the spectral condition.

Those remarks permit an elementary and constructive proof of the invertibility
of $\FE$. It suffices to show that the $\mathbb{C}[J]$-module $\mro{\rho}$
has rank $d$. That it cannot have rank greater than $d$ follows
quickly from the fact that a nonconstant function holomorphic on $\uhp/\SL$
must have poles at the cusps. It is enough then to find $d$ linearly
independent vectors in $\mro{\rho}$. Introduce the following notation:
given a $q$-series $f(q)=q^{\ell}\sum_{n=0}^{\infty}a_{n}q^{n}$
with $a_{0}\ne0$, define $o(f)$ to be $\ell$, the order of the
zero at $q=0$ -- e.g. $o(\eta)=1/24$ and $o(J)=-1$. The paper \cite{bib:KM}
explicitly constructs some weight $k$ vector-valued modular forms
for $\rho$, namely the Poincar\'{e} series $P$, where $k$ here
can be e.g. any sufficiently large multiple of 12. In particular,
let $\mathbb{Y}^{(i)}(\tau)=P(\tau;\rho,k,1,-2,i)$ in their notation,
for $1\le i\le d$; then each $\mathbb{Y}^{(i)}$ is a vector-valued
modular form for $\rho$ of weight $k$, holomorphic throughout $\uhp$,
with $o(\mathbb{Y}_{i}^{(i)})<0<o(\mathbb{Y}_{j}^{(i)})$ for all
$j\ne i$. Thus each $\mathbb{X}^{(i)}=\mathbb{Y}^{(i)}/\Delta^{k/12}$
lies in $\mro{\rho}$; that they are all linearly independent over
$\mathbb{C}(J)$ follows from the usual determinant argument.

\section{Low dimensional examples\label{sec:Low-dimensional-examples}}

This section is included to illustrate the effectiveness of the theory
on some simple examples up to dimension 3. As we shall see, some of
the nontrivial aspects of the theory already arise in these cases.
As usual, $S=\rho\!\sm{0}{-1}{1}{0}$ and $T=\rho\!\sm{1}{1}{0}{1}$
will denote the matrices representing the standard generators of $\SL$,
and $U=ST^{-1}=\rho\!\sm{0}{-1}{1}{-1}$.

The first comment is that
it is enough to consider indecomposable representations: 
indeed, if $\rho_{1}$ and $\rho_{2}$ are two representations of
$\SL$ satisfying the criteria of Section \ref{sec:The-fundamental-matrix},
then their direct sum $\rho_{1}\oplus\rho_{2}$ also satisfies these
criteria, and its exponent, characteristic and fundamental matrices
are just the direct sums of the corresponding matrices of its summands:\begin{align}
\FB\!\left(\rho_{1}\oplus\rho_{2}\right)= & \:\FB\!\left(\rho_{1}\right)\oplus\FB\!\left(\rho_{2}\right)\:,\nonumber \\
\mathcal{X}\!\left(\rho_{1}\oplus\rho_{2}\right)= & \:\mathcal{X\!}\left(\rho_{1}\right)\oplus\mathcal{X}\!\left(\rho_{2}\right)\:,\label{eq:Dsum}\\
\om{\rho_{1}\oplus\rho_{2}}= & \:\om{\rho_{1}}\oplus\om{\rho_{2}}\:.\nonumber \end{align}
 Thus, in order to determine the above quantities for an arbitrary
representation $\rho$, one should first decompose $\rho$ into a
direct sum of indecomposable representations, and determine the relevant
quantities for all the indecomposable constituents separately.

The representations of $\SL$ of dimension 1 that satisfy our criteria
are easy to classify: in this case the representation matrices are
mere numbers, and we get a total of six inequivalent representations,
each of which is a tensor power of the representation $\kan$ defined
in the last row of Table 1.
Note that this is in complete accord with the spectral condition:
there are exactly six pairs of 1-by-1 matrices that satisfy it. The
corresponding exponent and characteristic matrices are easily determined,
and this leads, via the compatibility equation Eq.\eqref{eq:compat},
to the corresponding fundamental matrices%
\footnote{There's no need to solve the differential equation in this case: the
fundamental matrices can be determined by purely function theoretic
arguments, or even better, from the determinantal formula Eq.\eqref{eq:detxi}.%
}. The results are gathered in Table 1, where $\omega=\exp\left(2\pi i/6\right)$
and $E_{k}$ stands for the Eisenstein series of weight $k$. %
\begin{table*}

\caption{One dimensional representations}

\begin{tabular}{|c|c|c|c|c|r|r|r|c|}
\hline 
&
&
&
&
&
&
&
&
\tabularnewline
$\am$&
$\bm$&
$\FB$&
$\mathcal{X}$&
$\om{\tau}$&
$S$&
$T$&
$U$&
name\tabularnewline
&
&
&
&
&
&
&
&
\tabularnewline
\hline
\hline 
&
&
&
&
&
&
&
&
\tabularnewline
0&
0&
1&
0&
1&
1&
1&
1&
1\tabularnewline
&
&
&
&
&
&
&
&
\tabularnewline
\hline 
&
&
&
&
&
&
&
&
\tabularnewline
0&
1&
2/3&
248&
${\displaystyle \frac{E_{4}}{\Delta^{1/3}}}=\left(J+744\right)^{1/3}$&
1&
$\omega^{4}$&
$\omega^{2}$&
$\kan^{2}$\tabularnewline
&
&
&
&
&
&
&
&
\tabularnewline
\hline 
&
&
&
&
&
&
&
&
\tabularnewline
0&
2&
1/3&
496&
${\displaystyle \frac{E_{8}}{\Delta^{2/3}}}=\left(J+744\right)^{2/3}$&
1&
$\omega^{2}$&
$\omega^{4}$&
$\kan^{4}$\tabularnewline
&
&
&
&
&
&
&
&
\tabularnewline
\hline 
&
&
&
&
&
&
&
&
\tabularnewline
1&
0&
1/2&
-492&
${\displaystyle \frac{E_{6}}{\Delta^{1/2}}}=\left(J-984\right)^{1/2}$&
-1&
-1&
1&
$\kan^{3}$\tabularnewline
&
&
&
&
&
&
&
&
\tabularnewline
\hline 
&
&
&
&
&
&
&
&
\tabularnewline
1&
1&
1/6&
-244&
${\displaystyle \frac{E_{10}}{\Delta^{5/6}}}=\left(J+744\right)^{1/3}\left(J-984\right)^{1/2}$&
-1&
$\omega$&
$\omega^{2}$&
$\overline{\kan}$\tabularnewline
&
&
&
&
&
&
&
&
\tabularnewline
\hline 
&
&
&
&
&
&
&
&
\tabularnewline
1&
2&
-1/6&
4&
${\displaystyle \frac{E_{14}}{\Delta^{7/6}}}=\left(J+744\right)^{2/3}\left(J-984\right)^{1/2}$&
-1&
$\omega^{5}$&
$\omega^{4}$&
$\kan$\tabularnewline
&
&
&
&
&
&
&
&
\tabularnewline
\hline
\end{tabular}
\end{table*}

The most interesting comments about Table 1 are related to the first
and last rows. In the first row we find the trivial representation,
and one would naively expect that the corresponding exponent matrix
is $0$. But this choice doesn't satisfy Eq.\eqref{eq:Riemann-Roch}:
we have to take $\FB=1$ according to our definitions. And indeed,
this choice is consistent with the fact that the constants belong
to $\mro{\rho}$ if $\rho$ is trivial. The last row is even more
interesting: naively, one would take $\FB=5/6$, but this would be
again in conflict with Eq.\eqref{eq:Riemann-Roch}; the correct value
is $\FB=-1/6$. Indeed, if one would have $\FB=5/6$, then $\Delta^{1/6}\SB{1}{1}$
would be a weight 2 modular form for the trivial representation, but
no such form exists, by classical arguments \cite{bib:Ap}.

Let's now turn to higher dimensions. The pairs of matrices \[
\begin{aligned}\FB=\frac{1}{24}\left(\begin{array}{rr}
17\\
 & 11\end{array}\right)\:,\quad & \mathcal{X}=\left(\begin{array}{rr}
133 & 1248\\
56 & -377\end{array}\right)\:\\
\FB=\frac{1}{24}\left(\begin{array}{rr}
23\\
 & 5\end{array}\right)\:,\quad & \mathcal{X}=\left(\begin{array}{rr}
3 & 26752\\
2 & -247\end{array}\right)\:\end{aligned}
\]
 both correspond to dimension 2 representations $\rho$ with the same
matrix \[
S={\displaystyle \frac{1}{\sqrt{2}}}\left(\begin{array}{rr}
1 & 1\\
1 & -1\end{array}\right)\]
 (of course, $T=\exp\left(2\pi i\FB\right)$ by definition). The corresponding
fundamental matrices have $q$-expansions \[
q^{\FB}\left(\begin{array}{cc}
q^{-1}+133+1673q+11914q^{2}+\ldots & 1248+49504q+806752q^{2}+\ldots\\
56+968q+7504q^{2}+\ldots & q^{-1}-377-22126q-422123q^{2}-\ldots\end{array}\right)\]
 and \[
q^{\FB}\left(\begin{array}{cc}
q^{-1}+3+4q+7q^{2}+\ldots & 26752+1734016q+46091264q^{2}+\ldots\\
2+2q+6q^{2}+\ldots & q^{-1}-247-86241q-4182736q^{2}-\ldots\end{array}\right)\:.\]
 They describe representations associated to the Wess-Zumino-Novikov-Witten
models \cite{bib:CFT} of level 1 based on the Lie algebras $E_{7}$
and $A_{1}$ respectively (whose dimensions 133 and 3 appear as $\mathcal{X}_{11}$,
and whose character vectors are given by the first columns of the
corresponding fundamental matrix).

More generally, the solution for an arbitrary two-dimensional $\SL$-represen\-ta\-tion
$\rho$ can be obtained in closed form -- for example, the fundamental
matrices can be expressed as linear combinations of classical hypergeometric
series.

In dimension 3, the sequence \[
\FB_{k}=\frac{1}{48}\left(\begin{array}{ccc}
47-2k\\
 & 23-2k\\
 &  & 2+4k\end{array}\right)\:,\]
 \[
\mathcal{X}_{k}=\left(\begin{array}{ccc}
k\left(2k+1\right) & \frac{1}{3}\left(31-2k\right)\left(9+2k\right)\left(25+2k\right) & 2^{12-k}\left(23-2k\right)\\
2k+1 & \left(11-k\right)\left(25+2k\right) & -2^{12-k}\\
2^{k} & -2^{k}\left(25+2k\right) & 2k-23\end{array}\right)\:,\]
 where $k$ is an integer in the range $0\leq k<12$, correspond to
representations that share the same matrix \[
S=\frac{1}{2}\left(\begin{array}{ccc}
1 & 1 & \sqrt{2}\\
1 & 1 & -\sqrt{2}\\
\sqrt{2} & -\sqrt{2} & 0\end{array}\right)\:.\]
 For $k=0$ one recovers the representation associated to the Ising
model \cite{bib:CFT}: in this case the fundamental matrix may be
expressed in terms of Weber functions \cite{bib:BG}.

What happens in these examples holds more generally: different $\SL$-repre\-sen\-ta\-tions
can have identical matrix $S$, but (in dimension $\le5$ \cite{bib:TW},
though not higher) an irreducible representation is completely determined
by $T$.

One striking feature of all the above examples is that their characteristic
matrices are all integral (in a suitable basis). This is far from
being trivial, since most solutions of Eq.(\ref{eq:monodromy}) have
irrational $\mathcal{X}$. Actually, the reason for using the pair
$\left(\FB,\mathcal{X}\right)$ to characterize the representation
$\rho$, instead of e.g. the pair $\left(\am,\bm\right)$, comes from
the observation that for representations $\rho$ which have a Conformal
Field Theory origin, there always seems to exist a basis in which
the characteristic matrix is integral. 
We'll explore this issue in Section \ref{sec:Positivity-and-integrality}.

\section{The inversion formula\label{sec:The-inversion-formula}}

We have seen above that the knowledge of the fundamental matrix $\om{\tau}$
allows the determination of all canonical basis vectors through solving
the recursion relations, and this in turn allows to determine the
unique element $\FA\in\mro{\rho}$ with a given principal part $\FE\FA$.
Actually, there exists an explicit inversion formula which gives $\FA$
in terms of $\FE\FA$ and the fundamental matrix.

\vspace{3mm}

\textbf{Inversion formula:} \textsl{for $\FA\!\left(q\right)\in\mro{\rho}$
with principal part $\FE\FA$}, \textsl{one has \begin{equation}
\FA\!\left(q\right)=\om{q}\frac{1}{2\pi i}\oint\frac{J^{\prime}\negmedspace\left(z\right)}{J\!\left(q\right)-J\!\left(z\right)}\om{z}^{-1}z^{\FB}\FE\FA\left(z\right)dz\:,\label{eq:inversion}\end{equation}
 where $J^{\prime}\!\!\left(z\right)=-z^{-2}+\sum_{n=1}^{\infty}nc\left(n\right)z^{n-1}$
is the derivative of $J$, and the integral is over a closed contour
encircling the origin and contained in the circle of radius $\left|q\right|$.}

\begin{proof}
Since the principal part map $\FE$ is linear, it is enough to prove
Eq.\eqref{eq:inversion} for the canonical basis vectors, in which
case it reads \begin{equation}
\left[\SB{\xi}{n}\left(q\right)\right]_{\eta}=\frac{1}{2\pi i}\oint\frac{z^{\FB_{\xi\xi}-n}J^{\prime}\!\!\left(z\right)}{J\!\left(q\right)-J\!\left(z\right)}\left[\om{q}\om{z}^{-1}\right]_{\eta\xi}dz\:.\label{eq:cbasis}\end{equation}
 To see that Eq.\eqref{eq:cbasis} holds, let's introduce the matrix
valued generating function \begin{equation}
\gfxi=\sum_{n=1}^{\infty}\left[\SB{\eta}{n}\left(q\right)\right]_{\xi}z^{n-1}\,\,.\label{eq:gfdef}\end{equation}
 As we'll see below, this series has a nonzero radius of convergence
around $z=0$, and thus defines a holomorphic function of $z$ in
a small enough neighborhood, for any fixed value of $q$. This means
that $z^{-n}\gfxi$ has a pole at $z=0$ whose residue is \begin{equation}
\left[\SB{\eta}{n}\left(q\right)\right]_{\xi}=\frac{1}{2\pi i}\oint z^{-n}\gfxi dz\:,\label{eq:resxi}\end{equation}
 by the residue theorem.

Multiplying both sides of the recursion relation Eq.(\ref{eq:recursion})
by $z^{m}$, and summing from $m=1$, one gets \begin{multline}
\gfxi-\om{q}_{\xi\eta}=zJ\!\left(q\right)\gfxi\\
-\sum_{m=1}^{\infty}\sum_{n=1}^{m-1}c\left(n\right)\left[\SB{\eta}{m-n}\left(q\right)\right]_{\xi}z^{m}-\sum_{\rho}\mathcal{X}_{\rho\eta}\left(z\right)\om{q}_{\xi\rho}\,\,\,,\label{eq:wrec2}\end{multline}
 where \begin{equation}
\mathcal{X}_{\xi\eta}\left(z\right)=\sum_{m=1}^{\infty}\xm{\xi}{\eta}{m}z^{m}\,\,\,.\label{eq:Xzdef}\end{equation}
 The double sum on the rhs. of Eq.\eqref{eq:wrec2} may be rearranged
as follows: \begin{multline*}
\sum_{m=1}^{\infty}\sum_{n=1}^{m-1}c\left(n\right)\left[\SB{\eta}{m-n}\left(q\right)\right]_{\xi}z^{m}=\sum_{n=1}^{\infty}\sum_{m=n+1}^{\infty}c\left(n\right)z^{n}\left[\SB{\eta}{m-n}\left(q\right)\right]_{\xi}z^{m-n}\\
=\sum_{n=1}^{\infty}\sum_{k=1}^{\infty}c\left(n\right)z^{n}\left[\SB{\eta}{k}\left(q\right)\right]_{\xi}z^{k}=z\left(J\!\left(z\right)-z^{-1}\right)\gfxi\,\,\,,\end{multline*}
 so that finally Eq.\eqref{eq:wrec2} reads \begin{equation}
z\left(J\!\left(q\right)-J\!\left(z\right)\right)\mathfrak{X}\left(q,z\right)=\om{q}\left(\mathcal{X}\!\left(z\right)-{\bf 1}\right)\,\,.\label{eq:wrec3}\end{equation}

We still have to determine the generating function $\mathcal{X}\!\left(z\right)$.
To do this, let's consider Eq.\eqref{eq:wrec3} in the limit when
$q$ approaches $z$: on the rhs. we get simply $\om{z}\left(\mathcal{X}\!\left(z\right)-{\bf 1}\right)$,
while on the lhs. all terms vanish because of the factor $\left(J\!\left(q\right)-J\!\left(z\right)\right)$,
except for those that are singular in $q$, which yield\begin{multline*}
\lim_{q\rightarrow z}\left\{ z\left(J\!\left(q\right)-J\!\left(z\right)\right)q^{\FB}\sum_{m=1}^{\infty}q^{-m}z^{m-1}\right\} =\\
z^{\FB+1}\lim_{q\rightarrow z}\left\{ \frac{J\!\left(q\right)-J\!\left(z\right)}{q-z}\right\} =z^{\FB+1}\D{J}\left(z\right)\,\,.\end{multline*}
 Note that the geometric sum is convergent for $\left|z\right|<\left|q\right|$.
All in all, we get \begin{equation}
\mathcal{X}\!\left(z\right)-{\bf 1}=\D{J}\left(z\right)\om{z}^{-1}z^{\mathbf{1}+\FB}\,\,\,.\label{eq:Xz}\end{equation}
 Inserting this last expression into Eq.\eqref{eq:wrec3}, we arrive
at \begin{equation}
\mathfrak{X}\left(q,z\right)=\frac{\D{J}\left(z\right)}{J\!\left(q\right)-J\!\left(z\right)}\om{q}\om{z}^{-1}z^{\FB}\:,\label{eq:wrec4}\end{equation}
 and this -- together with Eq.\eqref{eq:resxi} -- leads to the inversion
formula. Since the fundamental matrix is invertible except for the
elliptic points, Eq.\eqref{eq:wrec4} shows that the generating function
$\gfxi$ is indeed convergent for small enough $\left|z\right|<\left|q\right|$. 
\end{proof}
Let's stress that the above proof gives more than just the inversion
formula: it provides closed expressions for all the canonical basis
vectors, as well as for their generating function $\mathfrak{X}\left(q,z\right)$.
Incidentally, in the case of the trivial representation Eq.\eqref{eq:wrec4}
is related to the {}``bivarial transformation'' \cite{bib:No} of
Monstrous Moonshine.

\section{Positivity and integrality\label{sec:Positivity-and-integrality}}

The representations $\rho$ of most interest to us (coming from Conformal
Field Theories and Vertex Operator Algebras) have character vectors
$\FA\in\mro{\rho}$ which are dimensions of $\mathbb{Z}$-graded vector
spaces, and so their $q$-expansions Eq.\eqref{eq:qexp} have nonnegative
integer coefficients $\FA[n]$. In this section we find conditions
on $\rho$ for the existence of such $\FA$. Incidentally, this is
also why we choose $\FB$ and $\mathcal{X}$ for our fundamental data:
in the cases of most interest to us, $\mathcal{X}$ is integral.

Throughout this section, let $\rho$ be an indecomposable matrix representation
of P$\SL$, such that $T$ is diagonal and unitary. Call a nonzero
vector $\FA$ \textit{nonnegative} (resp. \textit{integral}) if all
its $q$-coefficients are nonnegative real numbers (resp. integral).
Recall the map $o(\sum_{n\ge0}a_{n}q^{n+\ell})=\ell$ of Section 3.
First, we give some easy conditions for nonnegativity.

\textbf{Nonnegativity test}: \textit{Suppose $\rho$ has a nonnegative
$\FA\in\mro{\rho}$. Then the matrix $S$ must have a strictly positive
eigenvector with eigenvalue 1. Suppose in addition there is a unique
component of $\FA$, call it $\FA_{0}(\tau)$, with a pole at $q=0$
of maximal order: i.e. $o(\FA_{0})<o(\FA_{\nu})$ for all $\nu\ne0$.
Then every entry in the 0-th column of $S$ must be a nonnegative
real number. }

This uniqueness assumption holds e.g. for any canonical basis vector;
it also holds for the character vector $\FA$ coming from a (unitary) Conformal
Field Theory, where it corresponds to the vacuum primary field.

The proof is easy. The eigenvector will be the vector $\FA(\tau)$
evaluated at $\tau=i$, i.e. $q=e^{-2\pi}$: it is positive because
$q>0$, and it has eigenvalue 1 because $\tau\mapsto-1/\tau$ fixes
$i$. Next, choose any $\eta$ such that $S_{\eta0}\ne0$; as $\tau$
approaches 0 along the imaginary axis, the component $\FA_{\eta}(\tau)$
remains manifestly positive. Applying $\tau\mapsto-1/\tau$, this
is equivalent to $\tau$ approaching $i\infty$ along the imaginary
axis (i.e. $q\rightarrow0$), of $\sum_{\mu}S_{\eta\mu}\FA_{\mu}(\tau)$.
But by the uniqueness hypothesis, this is dominated by the $\mu=0$
term. Hence positivity forces $S_{\eta0}\ge0$ for that $\eta$.

Most $\rho$ fail the first condition: e.g. measure-0 of 2-dimensional
and 4-dimensional representations, and 1/8-th of 3-dimensional ones,
satisfy it. The second condition is more powerful: e.g. it quickly
shows that any central charge $c<24$ Conformal Field Theories or
Vertex Operator Algebras with modular representation identical to
that of the Ising model, will have character vectors identical to
it. More generally, it implies that there will be only finitely many
possibilities for the character vectors of $c<24$ theories, with
fixed modular representation.

Now let's turn to integrality. As we will now explain, the existence
of integral $\FA$ leads us directly to representations $\rho$ whose
kernel is a congruence subgroup, i.e. ker$\,\rho$ contains some principal
congruence group \begin{equation}
\Gamma(N)=\{ A\in\SL\,|\, A\equiv1\ ({\rm mod}\ N)\}\ .\label{eq:prcong}\end{equation}
 Each component $\FA_{\eta}(\tau)$ of $\FA\in\mro{\rho}$ will be
a modular function for the kernel ker$\,\rho$, which we will require
here to be of finite index in $\SL$. Most such subgroups are noncongruence.
An example of a modular function for a noncongruence subgroup is \begin{equation}
\sqrt{{\frac{\eta(\tau)}{\eta(13\tau)}}}=q^{-{\frac{1}{4}}}(1-{\frac{1}{2}}q-{\frac{5}{8}}q^{2}-{\frac{5}{16}}q^{3}-{\frac{45}{128}}q^{4}+\cdots)\ .\end{equation}
 Although its Fourier coefficients are all rational, they have unbounded
denominator. Indeed, the following observation is due originally to
Atkin and Swinnerton-Dyer \cite{bib:AS}:

\textbf{Integrality conjecture:} \textit{Suppose $f(\tau)=q^{c}\sum_{n=0}^{\infty}a_{n}q^{n/b}$
is a modular function, holomorphic in $\uhp$, for some subgroup $G$
of $\SL$ with finite index, where $c$ is rational and $b$ is a
positive integer. If all coefficients $a_{n}$ are algebraic integers,
then $G$ is a congruence subgroup.}

Conversely, a modular function $f$ for $\Gamma(N)$ has a $q$-expansion
of the form \begin{equation}
f(\tau)=\sum_{n=-\infty}^{\infty}a_{n}q^{n/N}\ ,\label{eq:Nqexp}\end{equation}
 where $a_{n}=0$ for all but finitely many $n<0$; if $f$ is holomorphic
in $\uhp$, the denominators of its coefficients $a_{n}$ (if rational)
will be bounded. The integrality conjecture implies that $\rho$ can
have integral $\FA$ only if ker$\,\rho$ is congruence. That 
 the kernel is a congruence subgroup for a representation coming from
 Rational Conformal Field Theory was established in \cite{bib:Ba}.

Suppose for the remainder of this section that the kernel of $\rho$
contains some $\Gamma(N)$ -- in that case $N$ can be taken to be
the order of $T$. This implies that $\rho$ can equivalently be interpreted
as a representation of the finite group SL$_{2}(\mathbb{Z}_{N})$,
where $\mathbb{Z}_{N}=\mathbb{Z}/N\mathbb{Z}$. Incidentally, this
congruence subgroup hypothesis is straightforward to verify for any
given $\rho$, using the presentations of SL$_{2}(\mathbb{Z}[\frac{1}{p}])$
in \cite{bib:Hu}, but in practise a very convenient test is that
if ker$\,\rho$ is a congruence subgroup, then for all integers $\ell$
coprime to $N$, the diagonal entries of $T^{\ell^{2}}$ and $T$
are identical apart from order. To see this, let \begin{equation}
G_{\ell}=ST^{\frac{1}{\ell}}ST^{\ell}ST^{\frac{1}{\ell}}=\mr{\ell}{0}{0}{{\ell}^{-1}}\label{eq:Gell}\end{equation}
 where $\frac{1}{\ell}$ is the inverse of $\ell$ mod $N$; then
$G_{\ell}TG_{\ell}^{-1}=T^{\ell^{2}}$.

Now, any finite-dimensional representation of a finite group is equivalent
to one defined over some cyclotomic field $\mathbb{Q}_{L}=\mathbb{Q}[\xi_{L}]$,
where $\xi_{L}$ is the root of unity $e^{2\pi i/L}$. Replacing $N$
if necessary by multiple, we thus can (and will) assume that $\rho$
is a representation of SL$_{2}(\mathbb{Z}_{N})$, and all entries
of all matrices $\rho(\gamma)$ lie in $\mathbb{Q}_{N}$. Call any
such $\rho$ `$N$-defined'. Call $\FA$ \textit{rational} (resp.
$\mathbb{Q}_{N}$\textit{-rational}) if all coefficients in the $q$-expansions
of each component $\FA_{\eta}(\tau)$ are rational numbers (resp.
in $\mathbb{Q}_{N}$). It is known that if $\FA$ is rational and
ker$\,\rho$ is congruence, then some nonzero multiple $n\FA$ will
be integral. The remainder of this section is devoted to stating and
proving a necessary and sufficient condition for rationality. Not
surprisingly this involves the language of Galois.

For any $\ell$ coprime to $N$, let $\sigma_{\ell}\in\mathrm{Gal}(\mathbb{Q}_{N}/\mathbb{Q})$
be the Galois automorphism sending $\xi_{N}$ to $\xi_{N}^{\ell}$.
Let $\sigma_{\ell}$ act on any matrix $A\in M_{d\times d}(\mathbb{Q}_{N})$
entry-wise.

\textbf{Rationality test}: \textit{Let $\rho$ be $N$-defined, and
$\FA\in\mro{\rho}$ have components $\FA_{\eta}$ whose coefficients
$a_{\eta,n}$, $n\le0$, in Eq.\eqref{eq:Nqexp} are all rational.
Then $\FA$ is rational (hence a multiple is integral) iff for all
$\ell$ coprime to $N$, \begin{equation}
\sigma_{\ell}(S)=G_{\ell}S\ ,\label{eq:ratcond}\end{equation}
 where $G_{\ell}$ is defined in Eq.\eqref{eq:Gell}. In this case,
$G_{\ell}$ is a $\mathbb{Q}$-matrix and $S$ is real, and every
column of $\om{\tau}$ is rational.}

The starting point for proving this is the observation that any component
$\FA_{\eta}(\tau)$ of any vector $\FA\in\mro{\rho}$ is among other
things a modular function for $\Gamma(N)$. The theory of these functions
is quite rich (see e.g. Chapter 6 of \cite{bib:Sh} or Chapter 6 of
\cite{bib:La}).

Note that any $\FA\in\mro{\rho}$ is $\mathbb{Q}_{N}$-rational iff
all coefficients in the principal part are in $\mathbb{Q}_{N}$. In
particular, every canonical basis vector $\FA^{(\eta;m)}$ is $\mathbb{Q}_{N}$-rational.
The reason for this is that the space of modular forms for $\Gamma(N)$
of any weight $k$ has a basis with integral $q$-expansions, so so
does the space of modular functions for $\Gamma(N)$, holomorphic
in $\uhp$ and with bounded poles at the cusps; we can express $\FA_{\eta}(\tau)$
in terms of these basis functions by matching behaviours at the cusps,
and because $\rho$ is $N$-defined the coefficients will never leave
the field $\mathbb{Q}_{N}$.

The Galois automorphisms $\sigma_{\ell}$ mentioned above act on the
data $(\Lambda,\mathcal{X},\am,\bm,\om{\tau},\rho)$ associated to
any $N$-defined $\rho$, as follows. Note that the matrices in Eqs.(\ref{eq:ABdefa},b)
corresponding to $\Lambda$ and $\sigma_{\ell}\mathcal{X}$ will be
$\sigma_{\ell}\am$ and $\sigma_{\ell}\bm$, and thus the spectral
condition will be satisfied -- indeed the signature $(d,\alpha,\beta_{1},\beta_{2})$
won't have changed. It is easy to verify that the differential equation
Eq.\eqref{eq:hyper} will have solution $\sigma_{\ell}\om{\tau}$,
where we apply $\sigma_{\ell}$ entry-by-entry, and its action on
a $\mathbb{Q}_{N}$-rational $q$-series Eq.\eqref{eq:Nqexp} is simply
\begin{equation}
(\sigma_{\ell}f)(\tau)=\sum_{n=-\infty}^{\infty}\sigma_{\ell}(a_{n})q^{n/N}\ .\label{eq:galqser}\end{equation}
 By the above series, these $q$-series will be holomorphic throughout
$\uhp$. The corresponding P$\SL$-representation $\tilde{\rho}$
can be found by the following consideration.

Let $\mathcal{H}_{N}$ be the modular functions $f$ for $\Gamma(N)$,
holomorphic throughout $\uhp$, with coefficients $a_{n}\in\mathbb{Q}_{N}$.
The group GL$_{2}(\mathbb{Z}_{N})$ acts on $\mathcal{H}_{N}$ on
the right, i.e. $f|_{\alpha\circ\beta}=(f|_{\alpha})|_{\beta}$, as
follows (see Section 6.3 of \cite{bib:La} for more details). GL$_{2}(\mathbb{Z}_{N})$
is generated by SL$_{2}(\mathbb{Z}_{N})$, together with all matrices
of the form $M_{\ell}=\sm{1}{0}{0}{\ell}$ where $\ell$ is coprime
to $N$. $\gamma\in SL_{2}(\mathbb{Z}_{N})$ acts on $\mathcal{H}_{N}$
in the obvious way: first lift to $\SL$, then act on $\tau$ by that
fractional linear transformation. Moreover, $f|_{M_{\ell}}=\sigma_{\ell}f$,
as given by Eq.\eqref{eq:galqser}, recovering the action on $\om{\tau}$
we obtained last paragraph. That $\sigma_{\ell}f$ is holomorphic
in $\uhp$ iff $f$ is, follows from the previous paragraph (though
this is presumably also known classically). Then, writing $A=\sm{a}{b}{c}{d}$,
we have the calculation \begin{multline*}
(\sigma_{\ell}\FA)\left({\frac{a\tau+b}{c\tau+d}}\right)=\FA|_{M_{\ell}A}(\tau)=\sigma_{\ell}(\FA|_{M_{\ell}AM_{\ell}^{-1}})(\tau)\\
=\sigma_{\ell}(\mr{a}{\ell^{-1}b}{\ell c}{d}\FA(\tau))=(\sigma_{\ell}\mr{a}{\ell^{-1}b}{\ell c}{d})\sigma_{\ell}\FA(\tau)\end{multline*}
 where $\ell^{-1}$ denotes the inverse of $\ell$ mod $N$. Hence
we obtain\begin{equation}
\tilde{\rho}\sm{a}{b}{c}{d}=\sigma_{\ell}(\rho\sm{a}{\ell^{-1}b}{\ell c}{d})\ .\label{eq:rhotilde}\end{equation}
 More generally, if $\FA\in\mro{\rho}$ is $\mathbb{Q}_{N}$-rational,
then the same argument shows that $\sigma_{\ell}\FA$ lies in $\mro{\tilde{\rho}}$.

To complete the proof of the rationality test, note that for all $\ell$
coprime to $N$, $\sigma_{\ell}\FA\in\mro{\rho}$ iff \begin{equation}
\mr{a}{b}{c}{d}=\sigma_{\ell}(\mr{a}{\ell^{-1}b}{\ell c}{d})\ .\label{eq:****}\end{equation}
 Each component $(\sigma_{\ell}\FA)_{\eta}$ will have the same coefficients
$a_{n}$ as $\FA_{\eta}$, for all $n\le0$, and will be holomorphic
in $\uhp$. Hence $\sigma_{\ell}\FA=\FA$ for all $\ell$, i.e. $\FA$
is rational. Now, it suffices to test condition Eq.\eqref{eq:****}
at the generators $S$ and $T$. One leads to Eq.\eqref{eq:ratcond},
and the other to $T=\sigma_{\ell}T^{\ell^{-1}}$, which is automatically
satisfied. That $S$ is real follows from complex conjugation $\ell={-1}$
in Eq.\eqref{eq:ratcond}. That $G_{\ell}$ is rational follows from
the calculation \[
G_{\ell'}G_{\ell}S=\sigma_{\ell'\ell}S=(\sigma_{\ell'}G_{\ell})G_{\ell'}S=\sigma_{\ell'}(G_{\ell}G_{\ell'})S=\sigma_{\ell'}(G_{\ell'}G_{\ell})S\ .\]

The condition Eq.\eqref{eq:ratcond} is automatic in Conformal Field
Theory -- in this case $G_{\ell}$ is in fact monomial.

\section{Summary and outlook}

This paper solves the Riemann-Hilbert problem for P$\SL$: given a
representation $\rho$, we have a differential equation Eq.\eqref{eq:hyper}
whose monodromy is determined by $\rho$. The solution of this differential
equation  is the fundamental matrix $\om{\tau}$ of Eq.\eqref{eq:funddef}
-- given it, any vector-valued modular function $\FA$ with multiplier
$\rho$ can be uniquely determined from the inversion formula Eq.\eqref{eq:inversion}.
As an application of this theory, explicit bases for -- and dimensions
of -- spaces of vector-valued modular forms of half-integer weight
can be found. In practice, the most interesting vector-valued modular
functions have nonnegative integer $q$-expansions; the consequences
for $\rho$ of the existence of such vectors is worked out in Section
6.

A number of future developments are suggested by the analysis of this
paper. It is tempting to guess that the theory developed here can
be extended to other genus-0 discrete subgroups of PSL$_{2}(\mathbb{R})$.
There are 6486 such groups with the additional property that they
contain some $\Gamma(N)$ \cite{bib:Cum}: roughly a third of these
have only one cusp -- these may be the ones most accessible to our
methods.

Although vector-valued modular forms of half-integer weight can be
easily reduced to the modular functions studied here, the extension
to arbitrary weight will take more work. Such modular forms arise
naturally in Conformal Field Theory, and so this extension should
be pursued. Knopp and Mason \cite{bib:KM} have addressed questions
like the asymptotic growth of Fourier coefficients of these modular
forms, with methods apparently more effective when the weight is higher.
Our results would complement theirs: we would obtain bases and dimensions
 for any weight. We would suspect that uniform statements here for arbitrary 
weight will involve the braid group.

One could also speculate about the possibility of considering infinite
dimensional representations of $\SL$, which appear for instance in
quasi-rational Conformal Field Theory. In this case an indirect approach
could prove fruitful: first, solve Eq.\eqref{eq:monodromy} in an
arbitrary Banach algebra, then consider the solutions of the corresponding
differential equation Eq.\eqref{eq:hyper}; of course, all relevant
quantities that make sense will take their value in the given Banach
algebra. The technicalities involved are far from being clear.

Integrality and positivity, already touched upon in Section 6, lead
to many deep questions. For example, not all choices of $\FB$ compatible
with the trace formula Eq.\eqref{eq:Riemann-Roch} are equally good:
integrality, for instance, can be gained or lost by transformations
as in Eq.(\ref{eq:trafoL}), as Eq.(\ref{eq:trafoM}) shows. It would
be interesting to understand better how to choose the most suitable
$\FB$ in this respect.


\appendix

\section{The reduction of the modular representation}

To any Rational CFT is associated a finite dimensional representation
$\varrho$ of $\SL$, where in general $\varrho\sm{-1}{0}{0}{-1}$ is
not the identity, but a permutation matrix (charge conjugation). Nevertheless,
exploiting the fact that characters of charge conjugate primaries
are equal, one can associate to such a $\varrho$ a representation
$\rho$ for which $\rho\sm{-1}{0}{0}{-1}$ is the unit matrix, so
that the results of the paper may be applied. As far as conformal
characters are concerned, it is only $\rho$ that matters.

The procedure is as follows: let $T=\varrho\sm{1}{1}{0}{1}$ and $S=\varrho\sm{0}{-1}{1}{0}$
as usual. We know that $S^{2}$ is
a permutation matrix of order 2, representing charge conjugation.
An orbit $\eta$ of  charge conjugation has either length $\left|\eta\right|=1$,
or length $\left|\eta\right|=2$. 
For any such orbit $\eta$ we select a representative $\eta^{*}\in\eta$.

Define matrices $\mathcal{T}$ and $\mathcal{S}$, whose rows and
columns are indexed by these orbits $\eta$, via the rule \begin{align}
\mathcal{T}_{\xi\eta} & =\delta_{\xi\eta}T_{\eta^{*}\eta^{*}}\:,\nonumber \\
\mathcal{S}_{\xi\eta} & =\sum_{p\in\eta}S_{\xi^{*}p}\:.\label{eq:redrho}\end{align}
These matrices are well defined,
i.e. independent of the choice of the representatives $\xi^{*}\in\xi$
(since $S^2$ commutes with both $T$ and $S$), and they
determine a representation of $\SL$ which is trivial on the center:
this is the reduced representation $\rho$. Note that
all matrix elements of $\mathcal{S}$ are real numbers.

Some important properties of the modular representation $\varrho$
carry over to the reduced representation $\rho$ (e.g. the diagonality
of the Dehn-twist $T$), while others (like the symmetry and unitarity of $S$)
don't. The representation $\rho$ is equivalent to the largest subrepresentation
of $\varrho$ trivial on $\sm{-1}{0}{0}{-1}$.

We note that, while the reduction process results clearly in loss
of information, this loss is not as dramatic as one might expect:
for example, it is possible to reconstruct from the knowledge of $\rho$
the charge conjugation and the real part of $S$, as well as the full
matrix $T$.

\end{document}